# Research on Hopf Bifurcation and Stability of Heterogeneous Lorenz System with Single Time Delay

Zhu Erxi [1,2,3], Xu Min [3], Pi Dechang [1]

[1] College of Computer science and technology, Nanjing University of Aeronautics & Astronautics, Jiangsu Nanjing, 211106;
[2] Department of Computer Science and Technology, Suzhou Institute of Information Technology Jiangsu Su Zhou 215200;
[3] Department of Electronic and Information Engineering, Jiangsu Institute of Information Technology Jiangsu Wu Xi 214153.

Corresponding author: Zhu Erxi

Abstract: Time-delay chaotic systems refer to the hyperchaotic systems with multiple positive Lyapunov exponents. It is characterized by more complex dynamics and a wider range of applications as compared to those non-time-delay chaotic systems. In a three-dimensional general Lorenz chaotic system, time delays can be applied at different positions to build multiple heterogeneous Lorenz systems with a single time delay. Despite the same equilibrium point for multiple heterogeneous Lorenz systems with single time delay, their stability and Hopf bifurcation conditions are different due to the difference in time delay position. In this paper, the theory of nonlinear dynamics is applied to investigate the stability of the heterogeneous single-time-delay Lorenz system at the zero equilibrium point and the conditions required for the occurrence of Hopf bifurcation. First of all, the equilibrium point of each heterogeneous Lorenz system is calculated, so as to determine the condition that only zero equilibrium point exists. Then, an analysis is conducted on the distribution of the corresponding characteristic equation roots at the zero equilibrium point of the system to obtain the critical point of time delay at which the system is asymptotically stable at the zero equilibrium point and the Hopf bifurcation. Finally, mathematical software is applied to carry out simulation verification. Heterogeneous Lorenz systems with time delay have potential applications in secure communication and other fields.

# 1 Introduction

In respect of nonlinear dynamic systems, Hopf bifurcation and stability analysis of time-delay chaotic systems has attracted much attention for research. The distinctive characteristics of a time-delay chaotic system are detailed as follows. Firstly, the evolution of the system over time depends not only on the current state of the system but also on its past state. Secondly, it contains an infinite dimensional state space and exhibits extremely complex dynamic behaviors, which makes it different from the non-time-delay chaotic system. For these reasons, it has been widely applied in such fields as natural science, engineering technology, and social science [1-3]. As far engineering application is concerned, it is common for time-delay chaotic systems to cause system instability and bifurcation in various forms. Among them, Hopf bifurcation [4-5] is a commonplace and most discussed. For the generation mechanism of Hopf bifurcation, the condition is that the equilibrium point of the system switches stably with the change to a certain system parameter while the nonlinearity of the system restricts the disrupted divergent motion to a narrow range. Therefore, the precondition for the existence of Hopf bifurcation can be determined by analyzing the distribution of characteristic roots, which indicates that the existence of a certain system parameter value makes the system characteristic equation have negative real parts except for a pair of single conjugate pure imaginary roots. Moreover, the parameter value is taken as the Hopf bifurcation point when the corresponding characteristic root curve meets the transversal condition.

In recent years, there is still little attention paid to the research on Hopf bifurcation of the Lorenz system with time delay. Since 1963 when the meteorologist Lorenz proposed the first classic Lorenz system, researchers have put forward various heterogeneous Lorenz systems [8-11], such as Lü system [6] and Liu [7] system and conducted analysis of its chaotic mechanism for application in engineering settings. Using the first Lyapunov coefficient, Mello et al. analyzed the bifurcation characteristics of the three-dimensional Lorenz-like system [13]. Li et al. investigated the bifurcation characteristics of a novel Loren-like chaotic system at different equilibrium points [14]. Wang et al. demonstrated the fractional bifurcation of a five-dimensional Lorenz-like system [15]. Besides, the Routh-Hurwitz criterion and the high-dimensional Hopf bifurcation theory were applied to study the Hopf bifurcation characteristics of the three-dimensional autonomous Lorenz system [16]. In general, the time-delay chaotic system equation is linearized at the singularity to obtain the transcendental equation. In this way, the distribution of the roots of the transcendental equation is relied on to determine the Hopf bifurcation condition of the time-delay chaotic system. Through an in-depth discussion conducted by J.K. Hale [17], a theoretical foundation is laid for the study of Hopf bifurcation in time-delay chaotic systems. Professor Wei Junjie et al. [18] applied Rouche's theorem to provide the zero-point distribution theorem of exponential polynomials, which promoted the research on Hopf bifurcation theory. Extending and applying the canonical type theory to delay differential equations, T.Faria and Magalhães proposed a canonical type calculation method, which contributed significantly to the development of bifurcation theory [19-21]. At present, the bifurcation research on time-delay chaotic systems has been on the rise gradually. By introducing a generalized form of a time-delayed Lorenz system (the Lorenz system has (2n+1) dimensions), Mahmoud Gamal et al. analyzed not only the stability of trivial fixed points and non-trivial fixed points but also the conditions required for the occurrence of Hopf bifurcation [22]. Kun et al. adopted an improved method of undetermined coefficients to verify the homoclinic orbit of the Chen system with linear time-delay feedback, based on which the spiral involute projection method was proposed [23]. Lian et al. [24] conducted research on the Hopf bifurcation of Lorenz-like systems with time delay. Li et al. explored Hopf bifurcation of disturbed Lorenz-like systems with time delay [25].

According to the research and analysis of aforementioned literatures, different time delay positions have a more significant impact on the dynamic behavior of the system. For the Lorenz Chaos system, no one has studied the stability and bifurcation conditions of Lorenz system with time delay from the point of view of time delay. In this paper, time delays can be applied at different positions to build multiple heterogeneous Lorenz systems with a single time delay in a three-dimensional general Lorenz chaotic system. Despite the same equilibrium point for multiple heterogeneous Lorenz systems with single time delay, their stability and Hopf bifurcation conditions show difference due to the different time delay positions. The stability of Lorenz system with heterogeneous single delay at zero equilibrium point and the condition of Hopf bifurcation are studied by nonlinear dynamics theory. The simulation results are consistent with the theoretical analysis.

This paper is structured as follows. In Section one, a brief introduction is made of the time-delay chaotic system and its Hopf bifurcation. In Section two, a general heterogeneous single-time-delay Lorenz system model is proposed, and the condition of only zero equilibrium point is indicated.

Section 3 elaborates on the Hopf bifurcation and stability conditions of the three types of heterogeneous Lorenz systems with a single time delay. Besides, mathematical software is adopted to carry out simulation verification, which reveals that the conclusions drawn are consistent with the results of theoretical analysis. Finally, the conclusions are detailed in the concluding section.

## 2 General Lorenz System Model

Proposed by Lü et al. in 2002, the unified chaotic system connects the Lorenz system, Lü system, and Chen system. Its system model is expressed as

$$\begin{cases} \dot{x} = (25\alpha + 10)(y - x) \\ \dot{y} = (28 - 35\alpha)x - xz + (29\alpha - 1)y \\ \dot{z} = xy - (8 + \alpha)z/3 \end{cases} \quad (1)$$

Where, $x$, $y$, and $z$ represent state variables, while $\alpha \in [0,1]$ represents the system parameter. When $\alpha \in [0, 0.8)$, $\alpha \in (0.8, 1]$, and $\alpha = 0.8$, the system is classed as the generalized Lorenz system, the generalized Chen system, and the generalized Lü system, respectively. The unified chaotic system model demonstrates the basic structure of Lorenz using a single parameter. However, the number of its system parameters are too small, thus limiting the parameter range. Then, researchers proposed the corresponding bifurcation laws and stability conditions through continuous updates by forming many variants of Lorenz chaotic systems [26] (such as Lorenz-like systems [27]). Without any compromise on generality, a general Lorenz system is proposed in this paper, and an investigation is conducted into the bifurcation law of its heterogeneous single-time-delay chaotic system. The dynamic equation of the system is expressed as follows:

$$\begin{cases} \dot{x} = a(y - x) \\ \dot{y} = bx + dy - xz \\ \dot{z} = -cz + xy \end{cases} \quad (2)$$

Where, $a$, $b$, $c$, and $d$ are system parameters. Eq. (2) involves 7 terms, among which there are only 2 nonlinear terms. Compared with other chaotic or hyperchaotic systems, the structure of this system is simpler, thus making it easier to implement the circuit. Therefore, the system is applicable in such fields as secure communication. Fig. 1 shows the trend of changes in the phase diagram and state vector of a general Lorenz system over time when $a = 10$, $b = 28$, $c = 8/3$, $d = -1$.

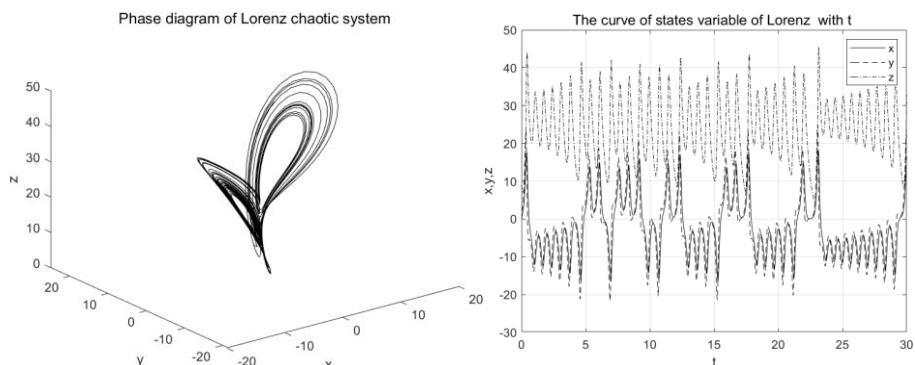

Fig. 1 The trend of changes in the general Lorenz system (the left exhibits the phase diagram of the Lorenz system in the $O - xyz$ space; the right indicates the change curve of $x$, $y$, $z$ with time $t$ for the Lorenz system)

As shown in Fig. 1, the general Lorenz system has two two-wing chaotic attractors when $a=10$, $b=28$, $c=8/3$, and $d=-1$. The state variable keeps oscillating, suggesting that the system has entered a state of chaos.

## 3 Stability and Hopf bifurcation conditions of heterogeneous single-time-delay Lorenz systems

Many researchers have imposed time delay on the state variables to develop the time-delay Lorenz chaotic system. However, the different positions of the added time delay can lead to a functional differential dynamic system with different dynamic behaviors. In Eq. (2), applying a single time delay to the state variable can give rise to nine forms of heterogeneous single-time-delay chaotic systems, which are different from other references. Three of the heterogeneous forms are presented as follows, which are also the stability and bifurcation conditions to be explored later in this study.

$$A.\begin{cases} \dot{x}=a(y(t-\tau)-x) \\ \dot{y}=bx+dy-xz \\ \dot{z}=-cz+xy \end{cases} \quad B.\begin{cases} \dot{x}=a(y-x(t-\tau)) \\ \dot{y}=bx+dy-xz \\ \dot{z}=-cz+xy \end{cases} \quad C.\begin{cases} \dot{x}=a(y-x) \\ \dot{y}=bx+dy-xz \\ \dot{z}=-cz(t-\tau)+xy \end{cases}$$

Where, $\tau(>0)$ denotes the amount of time delay, which can be understood as the time it takes for the predator to have the ability to prey, the incubation period of infectious diseases, or the delay time of signal transmission.

The three systems $ABC$ have three identical equilibrium points as follows:

$$(0,0,0), \quad (\sqrt{(b+d)c},\sqrt{(b+d)c},b+d), \quad (-\sqrt{(b+d)c},-\sqrt{(b+d)c},b+d)$$

When the three system parameters $a>0$, $b+d<0$, $c>0$, the system has a unique equilibrium point $O(0,0,0)$. The time-delay Lorenz system shows a better parameter range compared with other systems. Next, the stability of the equilibrium point $O(0,0,0)$ of three heterogeneous single-time-delay Lorenz systems is taken into consideration.

### 3.1 $A$ chaotic system

The linearization equation of the $A$ chaotic system at the equilibrium point $O(0,0,0)$ is expressed as:

$$\begin{cases} \dot{x}=a[y(t-\tau)-x] \\ \dot{y}=bx+dy \\ \dot{z}=-cz \end{cases} \quad (3)$$

Where, the value range of system parameters $a$, $b$, $c$, and $d$ is $a>0$, $b+d<0$, $c>0$, $a+c>d$, $|d|>|b|$. The corresponding characteristic equation of Eq. (3) is expressed as:

$$\begin{vmatrix} -a-\lambda & ae^{-\lambda\tau} & 0 \\ b & e-\lambda & 0 \\ 0 & 0 & -c-\lambda \end{vmatrix} = 0 \quad (4)$$

Eq. (4) can be reduced to

$$\lambda^3 + (a+c-d)\lambda^2 + (ca-cd-ad)\lambda - acd - ab(\lambda+c)e^{-\lambda\tau} = 0 \tag{5}$$

Where $p_1 = a+c-d$; $p_2 = ca-cd-ad$; $p_3 = -acd$; $p_4 = -ab$; $p_5 = -abc$. Thus, the following lemma can be obtained.

Lemma 1 If $\tau = 0$, the equilibrium point $O(0,0,0)$ of $A$ system is locally asymptotically stable when $a > 0$, $b+d < 0$, $c > 0$, $a+c > d$, and $|d| > |b|$.

Proof: when $\tau = 0$, the characteristic equation (5) is transformed into
$$\lambda^3 + p_1\lambda^2 + (p_2 + p_4)\lambda + p_3 + p_5 = 0 \tag{6}$$

Since $a > 0$, $b+d < 0$, $c > 0$, $a+c > d$, and $|d| > |b|$, it is easy to obtain that $p_1 > 0$, $p_2 + p_4 > 0$, and $p_3 + p_5 > 0$. According to the Routh-Hurwitz theorem, all roots of the characteristic equation (6) are common in having negative real parts. Thus, the equilibrium point $O(0,0,0)$ of $A$ system is asymptotically stable when $\tau = 0$.

When $\tau > 0$, suppose $\lambda = i\omega$ ($\omega$ is an undetermined constant greater than zero) is a pure imaginary root of Eq.(5), so that imaginary part $\omega$ satisfies
$$-i\omega^3 - p_1\omega^2 + ip_2\omega + p_3 + (ip_4\omega + p_5)(\cos\omega\tau - i\sin\omega\tau) = 0 \tag{7}$$

According to the equality of plural numbers, it can be obtained that:
$$\begin{cases} p_5 \cos\omega\tau + p_4\omega\sin\omega\tau = p_1\omega^2 - p_3 \\ p_4\omega\cos\omega\tau - p_5\sin\omega\tau = \omega^3 - p_2\omega \end{cases} \tag{8}$$

Eq. (8) can be equivalently transformed into
$$\omega^6 + (p_1^2 - 2p_2)\omega^4 + (p_2^2 - 2p_1p_3 - p_4^2)\omega^2 + p_3^2 - p_5^2 = 0 \tag{9}$$

A conclusion for Eq. (9) can be reached as follows.

Lemma 2 If $a > 0$, $b+d < 0$, $c > 0$, $a+c > d$, and $|d| > |b|$, Eq. (9) has at least one positive real root.

Proof: Set $u = \omega^2$, then, Eq. (9) can be reduced to
$$u^3 + (p_1^2 - 2p_2)u^2 + (p_2^2 - 2p_1p_3 - p_4^2)u + p_3^2 - p_5^2 = 0 \tag{10}$$

Suppose
$$f(u) = u^3 + (p_1^2 - 2p_2)u^2 + (p_2^2 - 2p_1p_3 - p_4^2)u + p_3^2 - p_5^2 \tag{11}$$

Eq. (11) can be converted into
$$f(u) = \frac{1 + (p_1^2 - 2p_2)\frac{1}{u} + (p_3^2 - 2p_1p_3 - p_4^2)\frac{1}{u^2} + (p_3^2 - p_5^2)\frac{1}{u^3}}{\frac{1}{u^3}} \tag{12}$$

It can be derived from Eq. (11) and Eq. (12) that
$$f(0) = p_3^2 - p_5^2 < 0, \quad \lim_{u \to +\infty} f(u) = +\infty$$

According to the theorem of the existence of function zeros, there is at least one real number $u_0 \in (0, +\infty)$ that makes $f(u_0) = 0$. Thus, Eq. (10) has one positive real root at minimum. Since $u = \omega^2$, Eq. (9) has at least one positive real root.

Suppose $\omega_0$ is a real root of Eq. (9), then Eq. (5) has a pure imaginary root $i\omega_0$. It can be obtained from Eq. (8) that

$$\cos\omega\tau = \frac{p_4\omega^4 - (p_1p_5 + p_2p_4)\omega^2 + p_3p_5}{p_4^2\omega^2 + p_5^2} \tag{13}$$

By substituting $\omega = \omega_0$ into Eq. (13), time delay $\tau$ can be calculated as

$$\tau_k = \frac{1}{\omega_0}\arccos(\frac{p_4\omega_0^4 - (p_2p_4 - p_1p_5)\omega_0^2 - p_3p_5}{p_4^2\omega_0^2 + p_5^2}) + \frac{2k\pi}{\omega_0}, k = 0,1,2,\cdots \quad (14)$$

Thus, $(\omega_0, \tau_k)$ is the solution of Eq. (5), suggesting that $\lambda = \pm i\omega_0$ is a pair of conjugate pure imaginary roots of Eq. (5) when $\tau = \tau_k$.

Suppose $\tau_0 = \min\{\tau_k\}$, then, time delay $\tau = \tau_0$ is the minimum value when the pure imaginary root $\lambda = \pm i\omega_0$ of Eq. (5) appears. Thus, there is a lemma shown as follows.

Lemma 3 If $a > 0$, $b + d < 0$, $c > 0$, $a + c > d$, $|d| > |b|$, and $\tau = \tau_0$, then, Eq. (5) has a pair of pure imaginary roots $\lambda = \pm i\omega_0$.

Suppose the characteristic root $\lambda(\tau) = \alpha(\tau) + i\omega(\tau)$ of Eq. (4) satisfies $\alpha(\tau_k) = 0$ and $\omega(\tau_k) = \omega_0$. The transversal conditions are presented below.

Lemma 4 If $a > 0$, $b + d < 0$, $c > 0$, $a + c > d$, $|d| > |b|$, and $f'(\omega_0^2) > 0$, then,

$$\frac{d\operatorname{Re}\lambda(\tau)}{d\tau}\bigg|_{\tau=\tau_k} > 0$$

Proof: The derivation regarding $\tau$ of both sides of Eq. (5) is performed to obtain

$$[3\lambda^2 + 2p_1\lambda + p_2 + p_4e^{-\lambda\tau} - \tau(p_4\lambda + p_5)e^{-\lambda\tau}]\frac{d\lambda}{d\tau} = \lambda(p_4\lambda + p_5)e^{-\lambda\tau} \quad (15)$$

It can be calculated according to Eq. (5) that

$$(p_4\lambda + p_5)e^{-\lambda\tau} = \lambda(\lambda^2 + p_1\lambda + p_2) \quad (16)$$

Substituting Eq. (16) into Eq. (15) yields

$$(\frac{d\lambda}{d\tau})^{-1} = -\frac{3\lambda^2 + 2p_1\lambda + p_2}{\lambda^2(\lambda^2 + p_1\lambda + p_2)} + \frac{p_4}{\lambda(p_4\lambda + p_5)} - \frac{\tau}{\lambda} \quad (17)$$

$\tau_k = i\omega_0$, therefore,

$$\operatorname{Re}[(\frac{d\lambda}{d\tau})^{-1}\bigg|_{\tau=\tau_k}] = -\operatorname{Re}[\frac{3\lambda^2 + 2p_1\lambda + p_2}{\lambda^2(\lambda^2 + p_1\lambda + p_2)}\bigg|_{\tau=\tau_k}] + \operatorname{Re}[\frac{p_4}{\lambda(p_4\lambda + p_5)}\bigg|_{\tau=\tau_k}]$$

$$= \operatorname{Re}[\frac{-3\omega_0^2 + 2ip_1\omega_0 + p_2}{\omega_0^2(\omega_0^2 - ip_1\omega_0 - p_2)}] + \operatorname{Re}(\frac{p_4}{p_4\omega_0^2 - ip_5\omega_0}) \quad (18)$$

$$= \frac{(p_2 - 3\omega_0^2)(\omega_0^2 - p_2) - 2p_1^2\omega_0^2}{\omega_0^2[(p_2 - \omega_0^2)^2 + p_1^2\omega_0^2]} - \frac{p_4^2}{p_4^2\omega_0^2 + p_5^2}$$

When $\tau = \tau_k$, Eq. (5) has pure imaginary roots $i\omega_0$, which are substituted into Eq. (5) to obtain

$$-i\omega_0^3 - p_1\omega_0^2 + ip_2\omega_0 + p_2 - (ip_4\omega_0 + p_5)e^{-i\omega_0\tau} = 0 \quad (19)$$

$|e^{-i\omega_0\tau}| = 1$ because $e^{-i\omega_0\tau} = \cos\omega_0\tau - i\sin\omega_0\tau$. Thus, it can be calculated using Eq. (19) that

$$|-p_1\omega_0^2 + p_3 + i(p_2\omega_0 - \omega_0^3)| = |-p_5 - ip_4\omega_0|$$

Namely,

$$\omega_0^2(p_2 - \omega_0^2)^2 + (p_1\omega_0^2 - p_3)^2 = (p_4\omega_0)^2 + p_5^2 \quad (20)$$

As obtained by combining Eq. (18) and Eq. (20),

$$\text{Re}[(\frac{d\lambda}{d\tau})^{-1}\big|_{\tau=\tau_k}] = \frac{3\omega_0^4 + 2(p_1^2 - 2p_2)\omega_0^2 + (p_2^2 - 2p_1p_3 - p_4^2)}{p_4^2\omega_0^2 + p_5^2} = \frac{f'(\omega_0^2)}{p_4^2\omega_0^2 + p_5^2} > 0$$

Besides, $Sign[\text{Re}(\frac{d\lambda}{d\tau}\big|_{\tau=\tau_k})] = Sign\{\text{Re}[(\frac{d\lambda}{d\tau})^{-1}\big|_{\tau=\tau_k}]\}$. Thus, the lemma is proved.

According to Lemma 4 and Hopf bifurcation theory, the following conclusions can be drawn.

Theorem 1 If $a > 0$, $b + d < 0$, $c > 0$, $a + c > d$, $|d| > |b|$, and $f'(\omega_0^2) > 0$, then,

(1) when $\tau \in [0, \tau_0)$, the equilibrium point $O(0,0,0)$ of $A$ system is asymptotically stable;

(2) when $\tau > \tau_0$, the equilibrium point $O(0,0,0)$ of $A$ system is unstable;

(3) $\tau = \tau_k (k = 0,1,2,\cdots)$ is the Hopf bifurcation value of $A$ system, suggesting that Hopf bifurcation occurs in $A$ system at the equilibrium point $O(0,0,0)$.

Considering that the parameters of $A$ system are $a > 0$, $b + d < 0$, $c > 0$, $a + c > d$, and $|d| > |b|$, $A$ system is simulated with $a = 10$, $b = -4$, $c = 2.5$, and $d = 2$. In this case, $A$ system can be converted into

$$\begin{cases} \dot{x} = 10y(t-\tau) - 10x \\ \dot{y} = -4x + 2y - xz \\ \dot{z} = -2.5z + xy \end{cases} \quad (21)$$

It can be calculated using mathematical software that the positive real root of Eq. (9) is $\omega_0 = 3.2376$, $f'(\omega_0^2) = 2.0909 \times 10^3 > 0$, and $\tau_0 = 0.2173$ in Eq. (14). Thus, Theorem 1 can be simplified into the following corollaries.

Corollary 1 If $a > 0$, $b + d < 0$, $c > 0$, $a + c > d$, $|d| > |b|$, and $f'(\omega_0^2) > 0$, then,

(1) when $\tau \in [0, 0.2173)$, the equilibrium point $O(0,0,0)$ of $A$ system is asymptotically stable;

(2) when $\tau > 0.2173$, the equilibrium point $O(0,0,0)$ of $A$ system is unstable;

(3) $\tau = 0.2173 + 0.6177k\pi (k = 0,1,2,3,\cdots)$ is the Hopf bifurcation value of $A$ system, suggesting that Hopf bifurcation occurs in $A$ system at the equilibrium point $O(0,0,0)$, leading to limit cycles.

Mathematical software is applied to draw the trajectory diagram and phase diagram of the state variable of $A$ system with time $t$ when the time delay $\tau$ takes different values, as illustrated in Fig. 2-4. The correctness of the results obtained is verified.

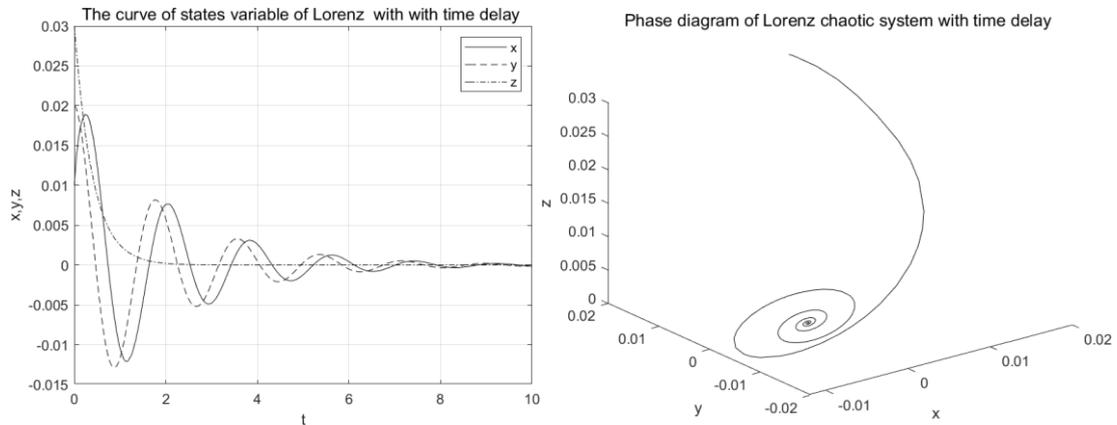

Fig. 2 The trend of changes in $A$ system when $\tau = 0.17$, $x(t) = 0.01$, $y(t) = 0.02$, and

$z(t) = 0.03 (t \in [-0.17, 0])$ (The left indicates the change curve of the state variables $x, y, z$ of A system with time t; the right exhibits the phase diagram of $A$ system in $O-xyz$ space)

As shown in Fig. 2, when $\tau = 0.17$, the value of the state variable $x, y, z$ of $A$ system approaches the equilibrium point $O(0,0,0)$ over time, as a result of which the equilibrium point $O(0,0,0)$ of $A$ system is asymptotically stable.

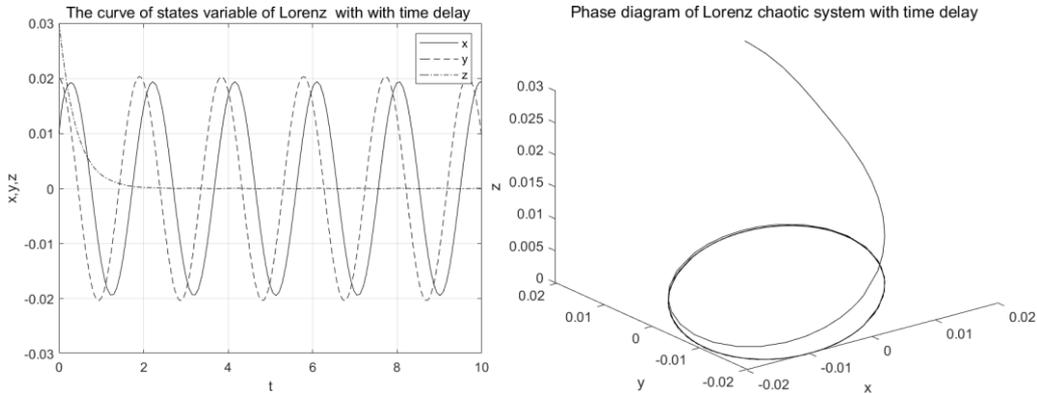

Fig. 3 The trend of changes in $A$ system when $\tau = 0.2173$, $x(t) = 0.01$, $y(t) = 0.02$, and $z(t) = 0.03 (t \in [-0.2173, 0])$ (The left indicates the change curve of the state variables $x, y, z$ of $A$ system with time t; the right presents the phase diagram of $A$ system in the $O-xyz$ space)

It can be observed in Fig. 3 that when $\tau = 0.2173$, the state variable $x, y, z$ of $A$ system keeps periodic oscillation with time $t$, and limit cycles appear in the $O-xyz$ space, suggesting that Hopf bifurcation occurs in $A$ system at the equilibrium point $O(0,0,0)$.

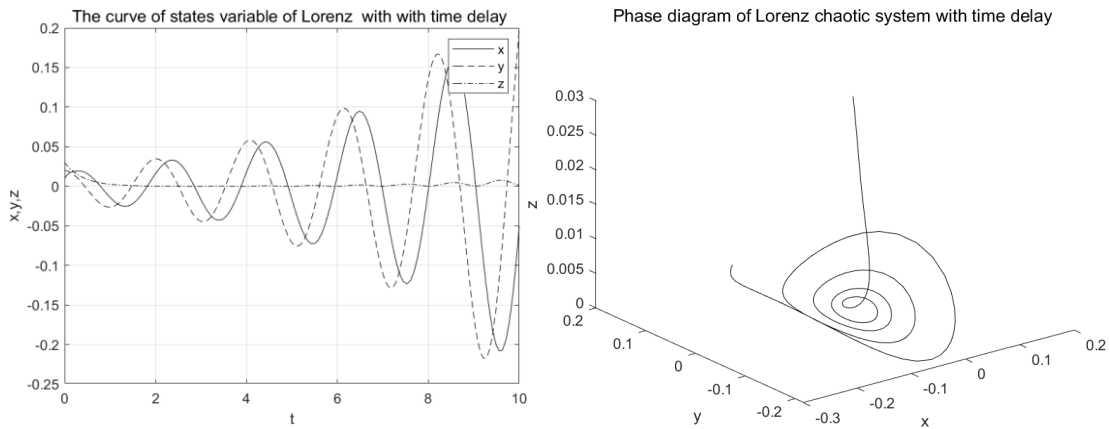

Fig. 4 The trend of changes in $A$ system when $\tau = 0.25$, $x(t) = 0.01$, $y(t) = 0.02$, and $z(t) = 0.03 (t \in [-0.25, 0])$. (The left indicates the change curve of the state variables $x, y, z$ of $A$ system with time t; the right presents the phase diagram of $B$ system in the $O-xyz$ space)

As shown from Fig. 4, the value of the state variables $x, y, z$ of $A$ system shifts away from the equilibrium point progressively with time $t$, suggesting that the equilibrium point $O(0,0,0)$ of $A$ system is unstable when $\tau = 0.25$.

## 3.2 B chaotic system

The linearization equation of the B chaotic system at the equilibrium point $O(0,0,0)$ is expressed as:

$$\begin{cases} \dot{x} = a(y - x(t-\tau)) \\ \dot{y} = bx + dy \\ \dot{z} = -cz \end{cases} \quad (22)$$

Where, the value range of system parameters $a$, $b$, $c$, and $d$ is $a > 0$, $b + d < 0$, $c > 0$, $a + c > d$, $|d| < |b|$. The corresponding characteristic equation of Eq. (22) is expressed as:

$$\begin{vmatrix} -ae^{-\lambda\tau} - \lambda & a & 0 \\ b & d - \lambda & 0 \\ 0 & 0 & -c - \lambda \end{vmatrix} = 0 \quad (23)$$

Eq. (23) can be reduced to

$$\lambda^3 + p_1\lambda^2 + p_2\lambda + p_5 + [a\lambda^2 + p_4\lambda + p_3]e^{-\lambda\tau} = 0 \quad (24)$$

Where $p_1 = c - d$; $p_2 = -(ab + cd)$; $p_3 = -acd$; $p_4 = a(c-d)$; $p_5 = -abc$. Thus, the following lemma can be obtained.

Lemma 5 If $\tau = 0$, the equilibrium point $O(0,0,0)$ of B system is locally asymptotically stable when $a > 0$, $b + d < 0$, $c > 0$, $a + c > d$, and $|d| < |b|$.

Proof: when $\tau = 0$, the characteristic equation (24) is transformed into

$$\lambda^3 + (a + p_1)\lambda^2 + (p_2 + p_4)\lambda + p_3 + p_5 = 0 \quad (25)$$

Since $a > 0$, $b + d < 0$, $c > 0$, $a + c > d$, and $|d| < |b|$, it is easy to obtain that $a + p_1 > 0$, $p_2 + p_4 > 0$, and $p_3 + p_5 > 0$. According to the Routh-Hurwitz theorem, all roots of the characteristic equation (6) are common in having negative real parts. Thus, the equilibrium point $O(0,0,0)$ of A system is asymptotically stable when $\tau = 0$.

When $\tau > 0$, suppose $\lambda = i\omega$ ($\omega$ is an undetermined constant greater than zero) is a pure imaginary root of Eq.(24), so that imaginary part $\omega$ satisfies

$$i\omega^3 + p_1\omega^2 - ip_2\omega - p_5 + (a\omega^2 - p_3 - ip_4\omega)(\cos\omega\tau - i\sin\omega\tau) = 0 \quad (26)$$

According to the equality of plural numbers, it can be obtained that:

$$\begin{cases} (a\omega^2 - p_3)\cos\omega\tau - p_4\omega\sin\omega\tau + p_1\omega^2 - p_5 = 0 \\ p_4\omega\cos\omega\tau + (a\omega^2 - p_3)\sin\omega\tau + p_2\omega - \omega^3 = 0 \end{cases} \quad (27)$$

Eq. (27) can be equivalently transformed into

$$\omega^6 + (p_1^2 - 2p_2 - a^2)\omega^4 + (p_2^2 - 2p_1p_5 + 2ap_3 - p_4^2)\omega^2 + p_5^2 - p_3^2 = 0 \quad (28)$$

A conclusion for Eq. (28) can be reached as follows.

Lemma 6 If $a > 0$, $b + d < 0$, $c > 0$, $a + c > d$, and $|d| < |b|$, Eq. (28) has at least one positive real root.

Proof: Set $u = \omega^2$, then, Eq. (28) can be reduced to

$$u^3 + (p_1^2 - 2p_2 - a^2)u^2 + (p_2^2 - 2p_1p_5 + 2ap_3 - p_4^2)u + p_5^2 - p_3^2 = 0 \quad (29)$$

Suppose

$$g(u) = u^3 + (p_1^2 - 2p_2 - a^2)u^2 + (p_2^2 - 2p_1p_5 + 2ap_3 - p_4^2)u + p_5^2 - p_3^2 \quad (30)$$

Eq. (30) can be converted into

$$g(u) = \frac{1 + (p_1^2 - 2p_2 - a^2)\frac{1}{u} + (p_2^2 - 2p_1p_5 + 2ap_3 - p_4^2)\frac{1}{u^2} + (p_5^2 - p_3^2)\frac{1}{u^3}}{\frac{1}{u^3}} \quad (31)$$

It can be derived from Eq. (30) and Eq. (31) that

$$g(0) = p_5^2 - p_3^2 < 0, \quad \lim_{u \to +\infty} g(u) = +\infty$$

According to the theorem of the existence of function zeros, there is at least one real number $u_0 \in (0, +\infty)$ that makes $g(u_0) = 0$. Thus, Eq. (29) has one positive real root at minimum. Since $u = \omega^2$, Eq. (28) has at least one positive real root.

Suppose $\omega_0$ is a real root of Eq. (29), then Eq. (28) has a pure imaginary root $i\omega_0$. It can be obtained from Eq. (27) that

$$\cos\omega\tau = \frac{(p_4 - ap_1)\omega^4 + (ap_5 + p_1p_3 - p_2p_4)\omega^2 - p_3p_5}{(a\omega^2 - p_3)^2 + p_4^2\omega^2} \quad (32)$$

By substituting $\omega = \omega_0$ into Eq. (32), time delay $\tau$ can be calculated as

$$\tau_k = \frac{1}{\omega_0}\arccos\left(\frac{m_1\omega_0^4 + m_2\omega_0^2 - p_3p_5}{(a\omega_0^2 - p_3)^2 + p_4^2\omega_0^2}\right) + \frac{2k\pi}{\omega_0}, k = 0, 1, 2, \cdots \quad (33)$$

where $m_1 = p_4 - ap_1$, $m_2 = ap_5 + p_1p_3 - p_2p_4$.

Thus, $(\omega_0, \tau_k)$ is the solution of Eq. (27), suggesting that $\lambda = \pm i\omega_0$ is a pair of conjugate pure imaginary roots of Eq. (24) when $\tau = \tau_k$.

Suppose $\tau_0 = \min\{\tau_k\}$, then, time delay $\tau = \tau_0$ is the minimum value when the pure imaginary root $\lambda = \pm i\omega_0$ of Eq. (24) appears. Thus, there is a lemma shown as follows.

Lemma 7 If $a > 0$, $b + d < 0$, $c > 0$, $a + c > d$, $|d| < |b|$, and $\tau = \tau_0$, then, Eq. (24) has a pair of pure imaginary roots $\lambda = \pm i\omega_0$.

Suppose the characteristic root $\lambda(\tau) = \alpha(\tau) + i\omega(\tau)$ of Eq. (24) satisfies $\alpha(\tau_k) = 0$ and $\omega(\tau_k) = \omega_0$. The transversal conditions are presented below.

Lemma 8 If $a > 0$, $b + d < 0$, $c > 0$, $a + c > d$, $|d| < |b|$, and $f'(\omega_0^2) > 0$, then,

$$\left.\frac{d\text{Re}\lambda(\tau)}{d\tau}\right|_{\tau=\tau_k} > 0$$

Proof: The derivation regarding $\tau$ of both sides of Eq. (24) is performed to obtain

$$[3\lambda^2 + 2p_1\lambda + p_2 + (2a\lambda + p_4)e^{-\lambda\tau} - \tau(a\lambda^2 + p_4\lambda + p_3)e^{-\lambda\tau}]\frac{d\lambda}{d\tau} \quad (34)$$
$$= \lambda(a\lambda^2 + p_4\lambda + p_3)e^{-\lambda\tau}$$

It can be calculated according to Eq. (24) that

$$(a\lambda^2 + p_4\lambda + p_3)e^{-\lambda\tau} = -(\lambda^3 + p_1\lambda^2 + p_2\lambda + p_5) \quad (35)$$

Substituting Eq. (35) into Eq. (34) yields

$$\left(\frac{d\lambda}{d\tau}\right)^{-1} = -\frac{3\lambda^2 + 2p_1\lambda + p_2}{\lambda(\lambda^3 + p_1\lambda^2 + p_2\lambda + p_5)} + \frac{2a\lambda + p_4}{\lambda(a\lambda^2 + p_4\lambda + p_3)} - \frac{\tau}{\lambda} \quad (36)$$

$\tau_k = i\omega_0$, therefore,

$$\mathrm{Re}[(\frac{d\lambda}{d\tau})^{-1}|_{\tau=\tau_k}] = -\mathrm{Re}[\frac{3\lambda^2 + 2p_1\lambda + p_2}{\lambda(\lambda^3 + p_1\lambda^2 + p_2\lambda + p_5)}|_{\tau=\tau_k}] + \mathrm{Re}[\frac{2a\lambda + p_4}{\lambda(a\lambda^2 + p_4\lambda + p_3)}|_{\tau=\tau_k}]$$

$$= -\mathrm{Re}[\frac{-3\omega_0^2 + 2ip_1\omega_0 + p_2}{\omega_0^4 - ip_1\omega_0^3 - p_2\omega_0^2 + ip_5\omega_0}] - \mathrm{Re}(\frac{i2a\omega_0 + p_4}{ia\omega_0^3 + p_4\omega_0^2 - ip_3\omega_0}) \quad (37)$$

$$= \frac{(3\omega_0^2 - p_2)(\omega_0^2 - p_2) + 2p_1(p_1\omega_0^2 - p_5)}{(p_5 - p_1\omega_0^2)^2 + [\omega_0^3 - p_2\omega_0]^2} - \frac{p_4^2 + 2a(a\omega_0^2 - p_3)}{p_4^2\omega_0^2 + (a\omega_0^2 - p_3)^2}$$

When $\tau = \tau_k$, Eq. (28) has pure imaginary roots $i\omega_0$, which are substituted into Eq. (24) to obtain

$$-i\omega_0^3 - p_1\omega_0^2 + ip_2\omega_0 + p_5 + (-a\omega_0^2 + ip_4\omega_0 + p_3)e^{-i\omega_0\tau} = 0 \quad (38)$$

$|e^{-i\omega_0\tau}| = 1$ because $e^{-i\omega_0\tau} = \cos\omega_0\tau - i\sin\omega_0\tau$. Thus, it can be calculated using Eq. (38) that

$$|-i\omega_0^3 - p_1\omega_0^2 + ip_2\omega_0 + p_5| = |-a\omega_0^2 + ip_4\omega_0 + p_3|$$

Namely,

$$(p_5 - p_1\omega_0^2)^2 + (p_2\omega_0 - \omega_0^3)^2 = (p_3 - a\omega_0^2)^2 + p_4^2\omega_0^2 \quad (39)$$

As obtained by combining Eq. (37) and Eq. (39),

$$\mathrm{Re}[(\frac{d\lambda}{d\tau})^{-1}|_{\tau=\tau_k}] = \frac{(3\omega_0^2 - p_2)(\omega_0^2 - p_2) + 2p_1(p_1\omega_0^2 - p_5) - p_4^2 - 2a(a\omega_0^2 - p_3)}{(p_5 - p_1\omega_0^2)^2 + (\omega_0^3 - p_2\omega_0)^2}$$

$$= \frac{f'(\omega_0^2)}{(p_5 - p_1\omega_0^2)^2 + (\omega_0^3 - p_2\omega_0)^2} > 0$$

Besides, $Sign[\mathrm{Re}(\frac{d\lambda}{d\tau}|_{\tau=\tau_k})] = Sign\{\mathrm{Re}[(\frac{d\lambda}{d\tau})^{-1}|_{\tau=\tau_k}]\}$. Thus, the lemma is proved.

According to Lemma 8 and Hopf bifurcation theory, the following conclusions can be drawn.

**Theorem 2** If $a > 0$, $b + d < 0$, $c > 0$, $a + c > d$, $|d| < |b|$, and $f'(\omega_0^2) > 0$, then,

(1) when $\tau \in [0, \tau_0)$, the equilibrium point $O(0,0,0)$ of $B$ system is asymptotically stable;

(2) when $\tau > \tau_0$, the equilibrium point $O(0,0,0)$ of $B$ system is unstable;

(3) $\tau = \tau_k (k = 0,1,2,\cdots)$ is the Hopf bifurcation value of $B$ system, suggesting that Hopf bifurcation occurs in $B$ system at the equilibrium point $O(0,0,0)$.

Considering that the parameters of $B$ system are $a > 0$, $b + d < 0$, $c > 0$, $a + c > d$, and $|d| < |b|$, $B$ system is simulated with $a = 10$, $b = 2$, $c = 2.5$, and $d = -4$. In this case, $B$ system can be converted into

$$\begin{cases} \dot{x} = 10y - 10x(t-\tau) \\ \dot{y} = 2x - 4y - xz \\ \dot{z} = -2.5z + xy \end{cases} \quad (40)$$

It can be calculated using mathematical software that the positive real root of Eq. (28) is $\omega_0 = 7.9396$, $f'(\omega_0^2) = 5.6868 \times 10^3 > 0$, and $\tau_0 = 0.18505$ in Eq. (28). Thus, Theorem 2 can be simplified into the following corollaries.

**Corollary 2** If $a > 0$, $b + d < 0$, $c > 0$, $a + c > d$, $|d| < |b|$, and $f'(\omega_0^2) > 0$, then,

(1) when $\tau \in [0, 0.18505)$, the equilibrium point $O(0,0,0)$ of $B$ system is asymptotically stable;

(2) when $\tau > 0.18505$, the equilibrium point $O(0,0,0)$ of $B$ system is unstable;

(3) $\tau = 0.18505 + 0.2519k\pi (k = 0,1,2,3,\cdots)$ is the Hopf bifurcation value of $B$ system, suggesting that Hopf bifurcation occurs in $B$ system at the equilibrium point $O(0,0,0)$, leading to limit cycles.

Mathematical software is applied to draw the trajectory diagram and phase diagram of the state variable of $B$ system with time $t$ when the time delay $\tau$ takes different values, as illustrated in Fig. 5-7. The correctness of the results obtained is verified.

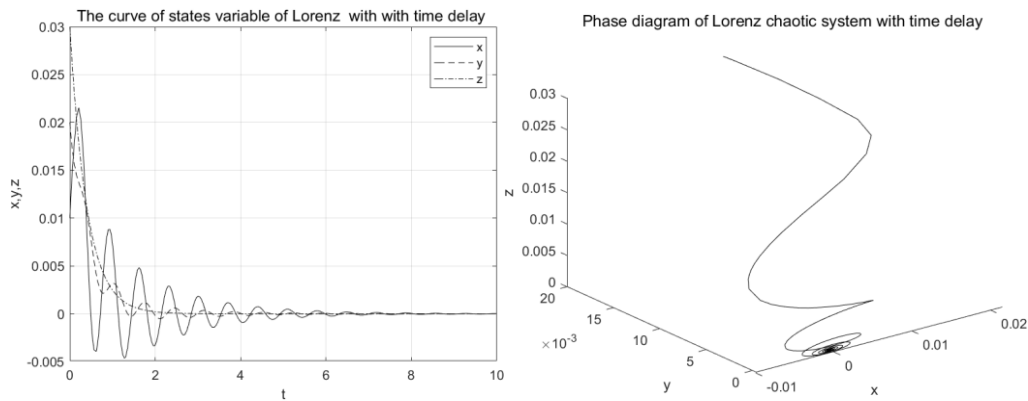

Fig. 5 The trend of changes in $B$ system when $\tau = 0.16$, $x(t) = 0.01$, $y(t) = 0.02$, and $z(t) = 0.03 (t \in [-0.16, 0])$ (The left indicates the change curve of the state variables $x, y, z$ of $B$ system with time t; the right exhibits the phase diagram of $B$ system in $O-xyz$ space)

As shown in Fig. 5, when $\tau = 0.16$, the value of the state variable $x, y, z$ of $B$ system approaches the equilibrium point $O(0,0,0)$ over time, as a result of which the equilibrium point $O(0,0,0)$ of $B$ system is asymptotically stable.

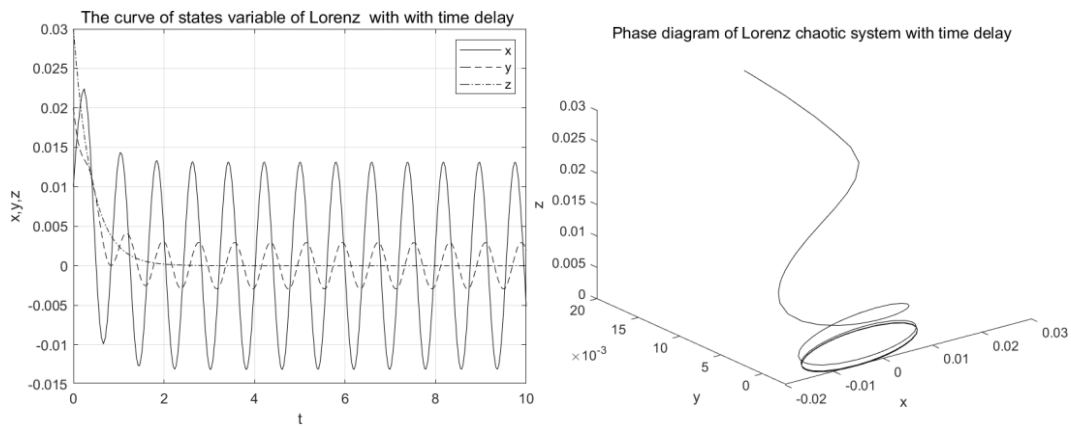

Fig. 6 The trend of changes in $B$ system when $\tau = 0.18505$, $x(t) = 0.01$, $y(t) = 0.02$, and $z(t) = 0.03 (t \in [-0.18505, 0])$ (The left indicates the change curve of the state variables $x, y, z$ of $B$ system with time t; the right presents the phase diagram of $B$ system in the $O-xyz$ space)

It can be observed in Fig. 6 that when $\tau = 0.18505$, the state variable $x, y, z$ of $B$ system keeps periodic oscillation with time $t$, and limit cycles appear in the $O-xyz$ space, suggesting that Hopf bifurcation occurs in $B$ system at the equilibrium point $O(0,0,0)$.

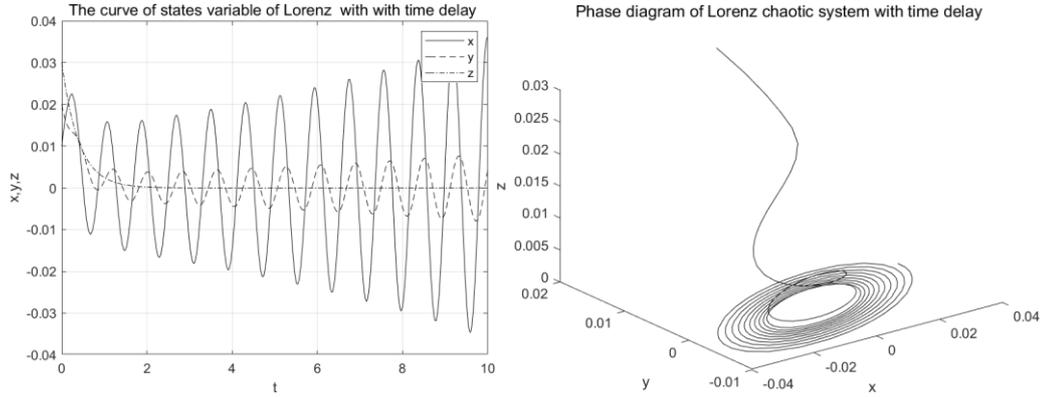

Fig. 7 The trend of changes in $B$ system when $\tau=0.19$, $x(t)=0.01$, $y(t)=0.02$, and $z(t)=0.03(t\in[-0.19,0])$ (The left indicates the change curve of the state variables $x, y, z$ of $B$ system with time t; the right presents the phase diagram of $B$ system in the $O-xyz$ space)

As shown from Fig. 7, the value of the state variables $x, y, z$ of $B$ system shifts away from the equilibrium point progressively with time $t$, suggesting that the equilibrium point $O(0,0,0)$ of $B$ system is unstable when $\tau=0.19$.

## 3.3 $C$ chaotic system

The linearization equation of the $C$ chaotic system at the equilibrium point $O(0,0,0)$ is expressed as:

$$\begin{cases} \dot{x}=a(y-x) \\ \dot{y}=bx+dy \\ \dot{z}=-cz(t-\tau) \end{cases} \tag{41}$$

Where, the value range of system parameters $a$, $b$, $c$, and $d$ is $a>0$, $b+d<0$, $c>0$, $a+c>d$. The corresponding characteristic equation of Eq. (41) is expressed as:

$$\begin{vmatrix} -a-\lambda & a & 0 \\ b & d-\lambda & 0 \\ 0 & 0 & -ce^{-\lambda\tau}-\lambda \end{vmatrix}=0 \tag{42}$$

Eq. (42) can be reduced to

$$\lambda^3+p_1\lambda^2+p_2\lambda+[c\lambda^2+p_3\lambda+p_4]e^{-\lambda\tau}=0 \tag{43}$$

Where $p_1=a-d$; $p_2=-a(b+d)$; $p_3=c(a-d)$; $p_4=-ac(b+d)$. Thus, the following lemma can be obtained.

Lemma 9 If $\tau=0$, the equilibrium point $O(0,0,0)$ of $A$ system is locally asymptotically stable when $a>0$, $b+d<0$, $c>0$ and $a+c>d$.

Proof: when $\tau=0$, the characteristic equation (43) is transformed into

$$\lambda^3+(c+p_1)\lambda^2+(p_3+p_2)\lambda+p_4=0 \tag{44}$$

Since $a>0$, $b+d<0$, $c>0$ and $a+c>d$, it is easy to obtain that $(c+p_1)>0$,

$p_3 + p_2 > 0$, and $p_4 > 0$. According to the Routh-Hurwitz theorem, all roots of the characteristic equation (44) are common in having negative real parts. Thus, the equilibrium point $O(0,0,0)$ of $C$ system is asymptotically stable when $\tau = 0$.

When $\tau > 0$, suppose $\lambda = i\omega$ ($\omega$ is an undetermined constant greater than zero) is a pure imaginary root of Eq. (43), so that imaginary part $\omega$ satisfies

$$i\omega^3 + p_1\omega^2 - ip_2\omega + [c\omega^2 - ip_3\omega - p_4](\cos\omega\tau - i\sin\omega\tau) = 0 \tag{45}$$

According to the equality of plural numbers, it can be obtained that:

$$\begin{cases} (p_4 - c\omega^2)\cos\omega\tau + p_3\omega\sin\omega\tau = p_1\omega^2 \\ p_3\omega\cos\omega\tau - (p_4 - c\omega^2)\sin\omega\tau = \omega^3 - p_2\omega \end{cases} \tag{46}$$

Eq. (46) can be equivalently transformed into

$$\omega^6 + (p_1^2 - 2p_2 - c^2)\omega^4 + (p_2^2 - 2cp_4 - p_3^2)\omega^2 - p_4^2 = 0 \tag{47}$$

A conclusion for Eq. (47) can be reached as follows.

Lemma 10 If $a > 0$, $b + d < 0$, $c > 0$ and $a + c > d$, Eq. (47) has at least one positive real root.

Proof: Set $u = \omega^2$, then, Eq. (47) can be reduced to

$$u^3 + (p_1^2 - 2p_2 - c^2)u^2 + (p_2^2 - 2cp_4 - p_3^2)u - p_4^2 = 0 \tag{48}$$

Suppose

$$h(u) = u^3 + (p_1^2 - 2p_2 - c^2)u^2 + (p_2^2 - 2cp_4 - p_3^2)u - p_4^2 \tag{49}$$

Eq. (49) can be converted into

$$h(u) = \frac{1 + (p_1^2 - 2p_2 - c^2)\frac{1}{u} + (p_2^2 - 2cp_4 - p_3^2)\frac{1}{u^2} - p_4^2\frac{1}{u^3}}{\frac{1}{u^3}} \tag{50}$$

It can be derived from Eq. (49) and Eq. (50) that

$$h(0) = -p_4^2 < 0, \quad \lim_{u \to +\infty} h(u) = +\infty$$

According to the theorem of the existence of function zeros, there is at least one real number $u_0 \in (0, +\infty)$ that makes $h(u_0) = 0$. Thus, Eq. (48) has one positive real root at minimum. Since $u = \omega^2$, Eq. (47) has at least one positive real root.

Suppose $\omega_0$ is a real root of Eq. (47), then Eq. (43) has a pure imaginary root $i\omega_0$. It can be obtained from Eq. (46) that

$$\cos\omega\tau = \frac{(p_4 - c\omega^2)p_1\omega^2 + p_3\omega(\omega^3 - p_2\omega)}{(p_4 - c\omega^2)^2 + p_3^2\omega^2} \tag{51}$$

By substituting $\omega = \omega_0$ into Eq.(51), time delay $\tau$ can be calculated as

$$\tau_k = \frac{1}{\omega_0}\arccos(\frac{(p_4 - c\omega_0^2)p_1\omega_0^2 + p_3\omega_0(\omega_0^3 - p_2\omega_0)}{(p_4 - c\omega_0^2)^2 + p_3^2\omega_0^2}) + \frac{2k\pi}{\omega_0}, k = 0, 1, 2, \cdots \tag{52}$$

Thus, $(\omega_0, \tau_k)$ is the solution of Eq. (43), suggesting that $\lambda = \pm i\omega_0$ is a pair of conjugate pure imaginary roots of Eq. (43) when $\tau = \tau_k$.

Suppose $\tau_0 = \min\{\tau_k\}$, then, time delay $\tau = \tau_0$ is the minimum value when the pure imaginary root $\lambda = \pm i\omega_0$ of Eq. (43) appears. Thus, there is a lemma shown as follows.

**Lemma 11** If $a > 0$, $b + d < 0$, $c > 0$, $a + c > d$ and $\tau = \tau_0$, then, Eq. (43) has a pair of pure imaginary roots $\lambda = \pm i\omega_0$.

Suppose the characteristic root $\lambda(\tau) = \alpha(\tau) + i\omega(\tau)$ of Eq. (43) satisfies $\alpha(\tau_k) = 0$ and $\omega(\tau_k) = \omega_0$. The transversal conditions are presented below.

**Lemma 12** If $a > 0$, $b + d < 0$, $c > 0$, $a + c > d$, and $f'(\omega_0^2) > 0$, then,

$$\frac{d \operatorname{Re} \lambda(\tau)}{d\tau}\bigg|_{\tau=\tau_k} > 0$$

Proof: The derivation regarding $\tau$ of both sides of Eq. (43) is performed to obtain

$$[3\lambda^2 + 2p_1\lambda + p_2 + (2c\lambda + p_3)e^{-\lambda\tau} - \tau(c\lambda^2 + p_3\lambda + p_4)e^{-\lambda\tau}]\frac{d\lambda}{d\tau} \quad (53)$$
$$= \lambda(c\lambda^2 + p_3\lambda + p_4)e^{-\lambda\tau}$$

It can be calculated according to Eq. (43) that

$$(c\lambda^2 + p_3\lambda + p_4)e^{-\lambda\tau} = -(\lambda^3 + p_1\lambda^2 + p_2\lambda) \quad (54)$$

Substituting Eq. (54) into Eq. (43) yields

$$(\frac{d\lambda}{d\tau})^{-1} = -\frac{3\lambda^2 + 2p_1\lambda + p_2}{\lambda(\lambda^3 + p_1\lambda^2 + p_2\lambda)} + \frac{2c\lambda + p_3}{\lambda(c\lambda^2 + p_3\lambda + p_4)} - \frac{\tau}{\lambda} \quad (55)$$

$\tau_k = i\omega_0$, therefore,

$$\operatorname{Re}[(\frac{d\lambda}{d\tau})^{-1}|_{\tau=\tau_k}] = -\operatorname{Re}[\frac{3\lambda^2 + 2p_1\lambda + p_2}{\lambda(\lambda^3 + p_1\lambda^2 + p_2\lambda)}|_{\tau=\tau_k}] + \operatorname{Re}[\frac{2c\lambda + p_3}{\lambda(c\lambda^2 + p_3\lambda + p_4)}|_{\tau=\tau_k}]$$
$$= -\operatorname{Re}[\frac{-3\omega_0^2 + 2ip_1\omega_0 + p_2}{\omega_0^4 - ip_1\omega_0^3 - p_2\omega_0^2}] + \operatorname{Re}(\frac{p_3 + i2c\omega_0}{-p_3\omega_0^2 - ic\omega_0^3 + ip_4\omega_0}) \quad (56)$$
$$= \frac{(3\omega_0^2 - p_2)(\omega_0^2 - p_2) + 2p_1^2\omega_0^2}{(\omega_0^3 - p_2\omega_0)^2 + p_1^2\omega_0^4} - \frac{p_3^2 + 2c(c\omega_0^2 - p_4)}{p_3^2\omega_0^2 + (c\omega_0^2 - p_4)^2}$$

When $\tau = \tau_k$, Eq. (43) has pure imaginary roots $i\omega_0$, which are substituted into Eq. (43) to obtain

$$-i\omega_0^3 - p_1\omega_0^2 + ip_2\omega_0 + (-c\omega_0^2 + ip_3\omega_0 + p_4)e^{-i\omega_0\tau} = 0 \quad (57)$$

$|e^{-i\omega_0\tau}| = 1$ because $e^{-i\omega_0\tau} = \cos\omega_0\tau - i\sin\omega_0\tau$. Thus, it can be calculated using Eq. (57) that

$$\left|-i\omega_0^3 - p_1\omega_0^2 + ip_2\omega_0\right| = \left|-c\omega_0^2 + ip_3\omega_0 + p_4\right|$$

Namely,

$$p_1^2\omega_0^4 + (\omega_0^3 - p_2\omega_0)^2 = (c\omega_0^2 - p_4)^2 + p_3^2\omega_0^2 \quad (58)$$

As obtained by combining Eq. (56) and Eq.(58),

$$\operatorname{Re}[(\frac{d\lambda}{d\tau})^{-1}|_{\tau=\tau_k}] = \frac{3\omega_0^4 + 2(p_1^2 - 2p_2 - c^2)\omega_0^2 + p_2^2 - 2cp_4 - p_3^2}{p_1^2\omega_0^4 + (\omega_0^3 - p_2\omega_0)^2}$$
$$= \frac{f'(\omega_0^2)}{p_1^2\omega_0^4 + (\omega_0^3 - p_2\omega_0)^2} > 0$$

Besides, $Sign[\text{Re}(\frac{d\lambda}{d\tau}|_{\tau=\tau_k})] = Sign\{\text{Re}[(\frac{d\lambda}{d\tau})^{-1}|_{\tau=\tau_k}]\}$. Thus, the lemma is proved.

According to Lemma 12 and Hopf bifurcation theory, the following conclusions can be drawn.

Theorem 1 If $a > 0$, $b+d < 0$, $c > 0$, $a+c > d$, and $h'(\omega_0^2) > 0$, then,

(1) when $\tau \in [0, \tau_0)$, the equilibrium point $O(0,0,0)$ of $C$ system is asymptotically stable;

(2) when $\tau > \tau_0$, the equilibrium point $O(0,0,0)$ of $C$ system is unstable;

(3) $\tau = \tau_k (k = 0,1,2,\cdots)$ is the Hopf bifurcation value of $A$ system, suggesting that Hopf bifurcation occurs in $C$ system at the equilibrium point $O(0,0,0)$.

Considering that the parameters of $C$ system are $a > 0$, $b+d < 0$, $c > 0$ and $a+c > d$, $C$ system is simulated with $a = 10$, $b = -4$, $c = 2.5$, and $d = 2$. In this case, $C$ system can be converted into

$$\begin{cases} \dot{x} = 10y - 10x \\ \dot{y} = -4x + 2y - xz \\ \dot{z} = -2.5z(t-\tau) + xy \end{cases} \quad (59)$$

It can be calculated using mathematical software that the positive real root of Eq. (47) is $\omega_0 = 3.7051$, $f'(\omega_0^2) = 8.0269 \times 10^2 > 0$, and $\tau_0 = 0.6265$ in Eq. (52). Thus, Theorem 3 can be simplified into the following corollaries.

Corollary 3 If $a > 0$, $b+d < 0$, $c > 0$, $a+c > d$, and $h'(\omega_0^2) > 0$, then,

(1) when $\tau \in [0, 0.6265)$, the equilibrium point $O(0,0,0)$ of $C$ system is asymptotically stable;

(2) when $\tau > 0.6265$, the equilibrium point $O(0,0,0)$ of $C$ system is unstable;

(3) $\tau = 0.6265 + 0.5398k\pi (k = 0,1,2,3,\cdots)$ is the Hopf bifurcation value of $C$ system, suggesting that Hopf bifurcation occurs in $C$ system at the equilibrium point $O(0,0,0)$, leading to limit cycles.

Mathematical software is applied to draw the trajectory diagram and phase diagram of the state variable of $C$ system with time $t$ when the time delay $\tau$ takes different values, as illustrated in Fig. 8-9. The correctness of the results obtained is verified.

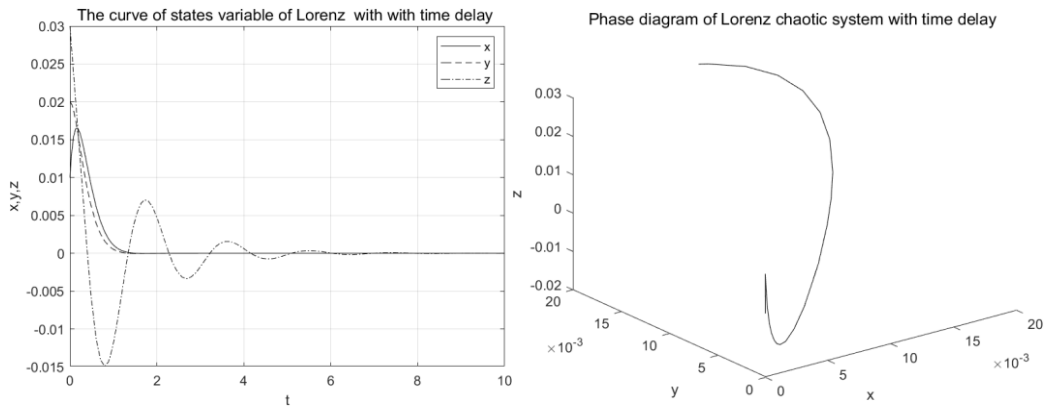

Fig. 8 The trend of changes in A system when $\tau = 0.4$, $x(t) = 0.01$, $y(t) = 0.02$, and $z(t) = 0.03 (t \in [-0.4, 0])$ (The left indicates the change curve of the state variables $x, y, z$ of $C$ system with time t; the right exhibits the phase diagram of A system in $O-xyz$ space)

As shown in Fig. 8, when $\tau = 0.4$, the value of the state variable $x, y, z$ of $C$ system approaches the equilibrium point $O(0,0,0)$ over time, as a result of which the equilibrium point $O(0,0,0)$ of $C$ system is asymptotically stable.

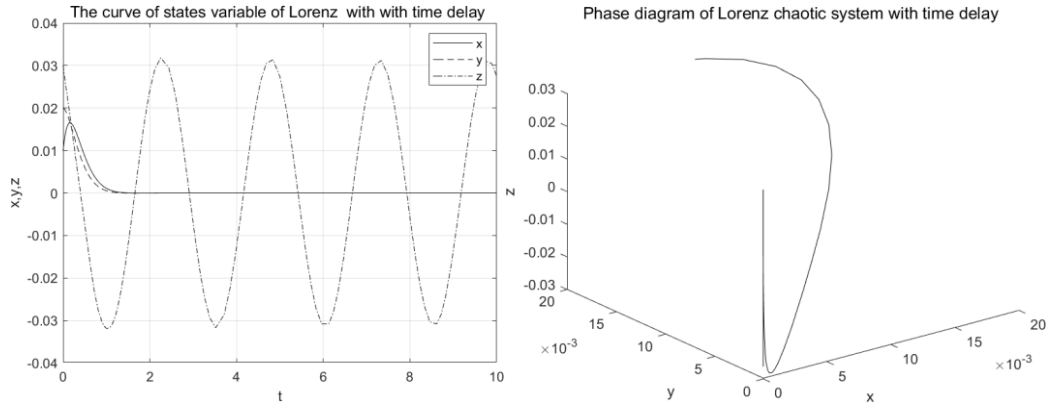

Fig. 9 The trend of changes in $C$ system when $\tau = 0.6265$, $x(t) = 0.01$, $y(t) = 0.02$, and $z(t) = 0.0001 (t \in [-0.6265, 0])$ (The left indicates the change curve of the state variables $x, y, z$ of A system with time t; the right presents the phase diagram of $C$ system in the $O - xyz$ space)

It can be observed in Fig. 9 that when $\tau = 0.6265$, the state variable $x, y, z$ of $C$ system keeps periodic oscillation with time $t$, and limit cycles appear in the $O - xyz$ space, suggesting that Hopf bifurcation occurs in $C$ system at the equilibrium point $O(0,0,0)$.

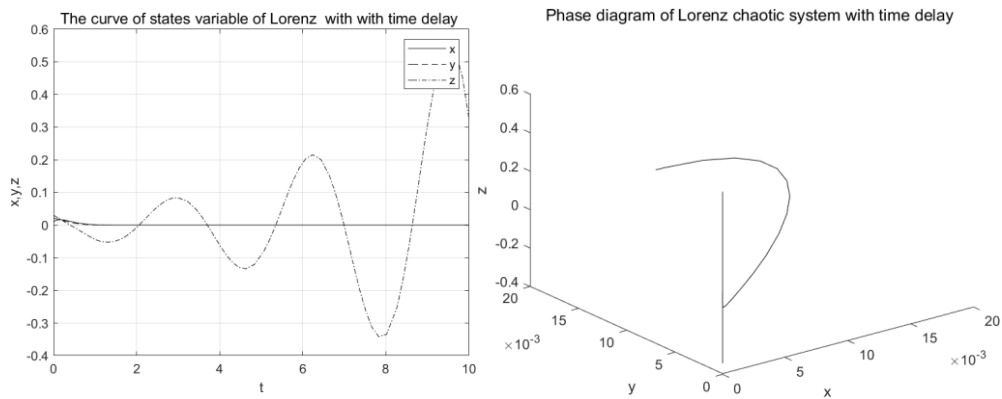

Fig. 10 The trend of changes in $C$ system when $\tau = 0.9$, $x(t) = 0.01$, $y(t) = 0.02$, and $z(t) = 0.03 (t \in [-0.9, 0])$. (The left indicates the change curve of the state variables $x, y, z$ of A system with time t; the right presents the phase diagram of $C$ system in the $O - xyz$ space)

In Fig. 10, the values of the state variables $x, y, z$ of the $C$ system gradually shift away from the equilibrium point with time $t$, suggesting that the equilibrium point $O(0,0,0)$ of $C$ system is unstable when $\tau = 0.9$. When $t < 32$, $x$ and $y$ state vectors keep approaching $O(0,0,0)$, and the amplitude of the state $z$ is on the increase. When $t \geq 32$, the state variables $x, y, z$ ceases to follow the original law and move into a state of chaos, as shown in Fig. 11.

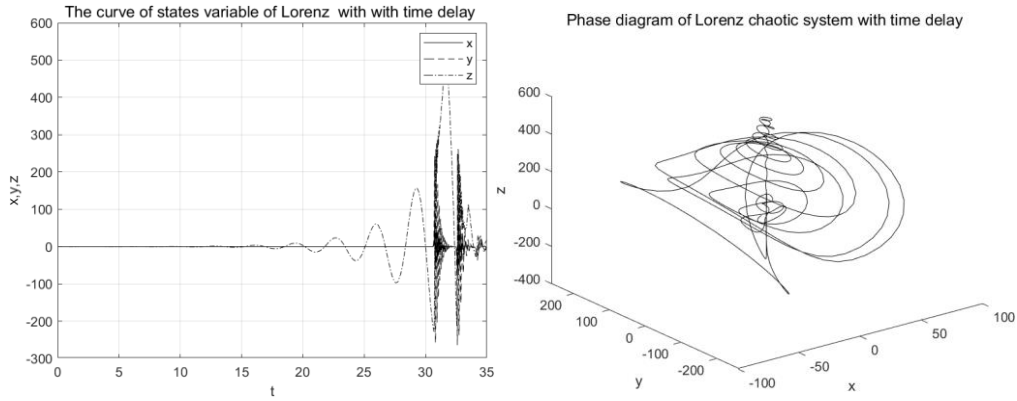

Fig. 11 The chaos state of C system when $\tau = 0.9$ and $t \geq 32$ (The left indicates the change curve of the state variables $x, y, z$ of C system with time t; the right shows the phase diagram of C system in the $O-xyz$ space.)

## 4. Conclusions

Since time delay is a common phenomenon in dynamic systems, it is necessary to explore the stability and bifurcation of functional differential dynamic systems. The position of time delay plays a significant role in the dynamic equation, with different time delay positions leading to different dynamic behaviors for the system. In this paper, a heterogeneous single-time-delay Lorenz system with different structures is constructed by loading time delays at different positions of the general Lorenz system. Three of the structures are selected to study the Hopf bifurcation and stability. According to the results, there is a single zero equilibrium point in the heterogeneous single-time-delay Lorenz system. Besides, the stability conditions at the zero equilibrium point and the parameter conditions required for the existence of Hopf bifurcation of different heterogeneous single-time-delay Lorenz systems are determined. Moreover, numerical simulation is performed to verify the correctness of the conclusions reached. The conclusion of this paper is expected to promote some existing literature research results.

## Conflicts of Interest

The author(s) declared no potential conflicts of interest with respect to the research, authorship, and/or publication of this article.

## Funding Statement

The author(s) received no financial support for the research, authorship, and/or publication of this article.

## Acknowledgments

Some of the authors of this publication are also working on these related projects: (1) higher vocational education teaching fusion production integration platform construction projects of Jiangsu province under grant no.2019(26), (2) natural science fund of Jiangsu province under grant

no.BK20131097,(3)"Qin Lan project" teaching team in colleges and universities of Jiangsu province under grant no.2017(15), (4)high level of Jiangsu province key construction project funding under grant no.2017(17).

## ORCID ID

Zhu Er Xi, https://orcid.org/0000-0003-1888-3377